\documentclass[letterpaper]{article}

\date{\today}

\newcommand{\myauthor}{Benjamin Antieau}
\newcommand{\mytitle}{On a theorem of Hazrat and Hoobler}

\title{\mytitle\footnote{This material is based upon work supported by the NSF under Grant RTG DMS 0838697.}}
\author{\myauthor}

\usepackage[pdfstartview=FitH,
            pdfauthor={\myauthor},
            pdftitle={\mytitle},
            colorlinks,
            linkcolor=reference,
            citecolor=citation,
            urlcolor=e-mail,
            backref]{hyperref}

\usepackage{rnote}
\usepackage{theorems}
\usepackage{times}
\usepackage[all,cmtip]{xy}
\usepackage{tikz}
\usetikzlibrary{matrix,arrows}

\let\oldmarginpar\marginpar
\renewcommand\marginpar[1]{\-\oldmarginpar[\raggedleft\footnotesize #1]%
{\raggedright\footnotesize #1}}

\begin{document}

\maketitle
\begin{abstract}
  \noindent
  We use cycle complexes with coefficients in an Azumaya algebra, as developed by Kahn and Levine, to compare the $G$-theory
  of an Azumaya algebra to the $G$-theory of the base scheme. We obtain a sharper version of a theorem of Hazrat and Hoobler in certain cases.

\paragraph{Mathematics Subject Classification 2000}
Primary: \href{http://www.ams.org/mathscinet/msc/msc2010.html?t=14Fxx&btn=Current}{14F22},
    \href{http://www.ams.org/mathscinet/msc/msc2010.html?t=19Dxx&btn=Current}{19Dxx}.
\end{abstract}

\section{Introduction}

Let $\K_*(X;\mathcal{A})$ be the $K$-theory of left $\mathcal{A}$-modules which are locally free and finite rank coherent $\mathcal{O}_X$-modules; let $\G_*(X;\mathcal{A})$
be the $K$-theory of left $\mathcal{A}$-modules which are coherent $\mathcal{O}_X$-modules.

We prove the following theorem.

\begin{theorem}
    Let $X$ be a $d$-dimensional scheme of finite type over a field $k$, and let $\mathcal{A}$ be an Azumaya algebra on $X$ of constant degree $n$.
    Let $B_{\mathcal{A}}:\G_i(X)\rightarrow\G_i(X;\mathcal{A})$ and
    $B_{\mathcal{A}}:\K_i(X)\rightarrow\K_i(X;\mathcal{A})$ be the homomorphisms induced by the functor $\mathcal{F}\mapsto \mathcal{A}\otimes_{\mathcal{O}_X}\mathcal{F}$. Then,
    \begin{enumerate}
        \item   the kernel and cokernel of $B_{\mathcal{A}}:\G_i(X)\rightarrow\G_i(X;\mathcal{A})$ are torsion groups of exponents dividing $n^{2d+2}$;
        \item   the kernel and cokernel of $B_{\mathcal{A}}:\K_i(X)\rightarrow\K_i(X;\mathcal{A})$ are torsion groups of exponents dividing $n^{2d+2}$ if $X$ is regular.
    \end{enumerate}
\end{theorem}

\begin{corollary}
    If $\mathcal{A}$ is an Azumaya algebra of constant degree $n$ over a scheme $X$ of finite type over a field $k$, then the base extension homomorphism
    \begin{equation*}
        B_{\mathcal{A}}:\G_*(X)\otimes_\ZZ \ZZ\left[\frac{1}{n}\right]\rightarrow\G_*(X;\mathcal{A})\otimes_\ZZ \ZZ\left[\frac{1}{n}\right]
    \end{equation*}
    is an isomorphism
\end{corollary}

The theorem above should be compared to the following two theorems, which motivated us in the first place.

\begin{theorem}[Hazrat-Millar~\cite{hazrat_millar}]
    If $\mathcal{A}$ is an Azumaya algebra of constant degree $n$ which is free over a noetherian affine scheme $X$, then
    \begin{equation*}
        B_{\mathcal{A}}:\K_i(X)\rightarrow\K_i(X;\mathcal{A})
    \end{equation*}
    has torsion kernel and cokernel of exponents at most $n^4$.
\end{theorem}

\begin{theorem}[Hazrat-Hoobler~\cite{hazrat_hoobler}]
    Let $X$ be a $d$-dimensional noetherian scheme, and let $\mathcal{A}$ be an Azumaya algebra on $X$ of constant degree $n$. Then,
    \begin{enumerate}
        \item   the kernel of $B_{\mathcal{A}}:\G_i(X)\rightarrow\G_i(X;\mathcal{A})$ is torsion of exponent dividing $n^{2d(d+1)+2}$, and the cokernel is torsion of exponent dividing $n^{4d+2}$;
        \item   the kernel of $B_{\mathcal{A}}:\K_i(X)\rightarrow\K_i(X;\mathcal{A})$ is torsion of exponent dividing $n^{2d(d+1)+2}$ if $X$ is regular,
            and the cokernel is torsion of exponent dividing $n^{4d+2}$ in this case;
        \item   the kernel and cokernel of $B_{\mathcal{A}}:\K_i(X)\rightarrow\K_i(X;\mathcal{A})$ are torsion groups of exponent dividing $n^{2d+2}$ if $X$ has an ample line bundle.
    \end{enumerate}
\end{theorem}

Since a degree $n$ Azumaya algebra is locally split by degree $n$ extensions, it is expected that the base extension map
\begin{equation}\label{eq:ghr}
    B_{\mathcal{A}}:\K_*(X)\otimes_\ZZ \ZZ\left[\frac{1}{n}\right]\rightarrow \K_*(X;\mathcal{A})\otimes_\ZZ \ZZ\left[\frac{1}{n}\right]
\end{equation}
should be an isomorphism.

Here is a partial history of results and techniques in this direction.

Wedderburn's theorem~\cite{herstein_noncommutative} easily implies that $\K_0(k)\rightarrow\K_0(A)$ is injective with cokernel isomorphic to $\ZZ/m$, where $A\iso M_m(D)$ for
a central $k$-division algebra $D$.

Green-Handelman-Roberts~\cite{green_handelman_roberts} proved that the map $B_{\mathcal{A}}$ is an isomorphism when $\mathcal{A}$ is a central simple algebra of degree $n$ over a field.
They used the Skolem-Noether theorem. That case has also been proven by Hazrat~\cite{hazrat_reduced} using the fact that $A$ is \'etale locally a matrix algebra.

The theorem of Hazrat-Millar quoted above uses the opposite algebra.
The theorem of Hazrat-Hoobler uses Bass-style stable range arguments and Zariksi descent for $G$-theory.

Our result uses twisted versions of Bloch's cycle complexes.
These twisted cycle complexes and the twisted motivic spectral sequence that relates them to $G$-theory are due to Kahn and Levine~\cite{kahn_levine}.
It is possible that our result could be extended to essentially smooth schemes over Dedekind rings by a combination of the work of Kahn and Levine~\cite{kahn_levine} and Geisser~\cite{geisser_dedekind}.

The following is an interesting corollary of our approach: there are natural filtrations of length $d$ on $\G_i(X)$
and $\G_i(X;\mathcal{A})$ coming from~\cite{kahn_levine}. The map $B_{\mathcal{A}}:\G_i(X)\rightarrow \G_i(X;\mathcal{A})$ respects the filtrations.
We show that the induced maps on each of the $d+1$ slices have kernel and cokernel groups of exponent at most $n^2$.

It is worth mentioning two other functors on Azumaya algebras with values in abelian groups where the base extension maps are isomorphisms.
Dwyer and Friedlander~\cite[2.4, 3.1]{dwyer_friedlander_az} showed that
\begin{equation*}
    \K_*^{\et}(R;\ZZ/m)\rightarrow\K_*^{\et}(R;A;\ZZ/m)
\end{equation*}
is an isomorphism in some cases (all of which are Azumaya algebras over a noetherian ring), where $\K^{\et}$ denotes \'etale $K$-theory, as, for instance, in Thomason~\cite{thomason_algebraic_1985}.
In this direction, it is possible to show (for instance, in the setting of Antieau~\cite{antieau_index_2009})
that $\K^{\et}(X;\mathcal{A})$ is an invertible object (in the sense of the Picard group) over $\K^{\et}(X)$ in the category of \'etale sheaves of $\K^{\et}$-module spectra
on a scheme $X$.

Finally, Corti\~nas and Weibel~\cite{cortinas_weibel} proved that the base extension maps induce isomorphisms in Hochschild homology over a field $k$.

\paragraph{Acknowledgments}
We thank Christian Haesemeyer, Roozbeh Hazrat, and Ray Hoobler for conversations.

\section{Twisted higher Chow groups and twisted $G$-theory}

Let $X$ in $\mathbf{Sch}/k$ be an integral $k$-scheme of finite type, and let $\mathcal{A}$ be a sheaf of
Azumaya algebras on $X$ of rank $n^2$. The degree of $\mathcal{A}$ is defined to be the integer $n$.
Let $\mathcal{E}$ be a left $\mathcal{A}$-module which is locally free and finite rank $na$ as an $\mathcal{O}_X$-module.
For generalities on Azumaya algebras, which as $\mathcal{O}_X$-modules are always locally free and of finite rank,
see~\cite{grothendieck_brauer_1}.

As in Kahn-Levine~\cite{kahn_levine}, define the cycle complex of
$X$ with coefficients in $\mathcal{A}$ as follows. Let $S^X_{(s)}(t)$ denote the set of closed
subsets $W\subset X\times_k\Delta^t$ such that
\begin{equation*}
    \dim_k W\cap X\times_k F\leq s+\dim_k F
\end{equation*}
for all faces $F$ of $\Delta^n$. Taking inverse images, $S^X_{(s)}(*)$ becomes a simplicial set.
Let $X_s(t)$ denote the subset of irreducible $W$ in $S^X_{(s)}(t)$ such that $\dim_k W=s+t$. Define, for $t\geq 0$,
\begin{equation*}
    z_s(X,t;\mathcal{A})=\bigoplus_{W\in X_s(t)}\K_0(k(W);\mathcal{A}).
\end{equation*}
See~\cite[Definition~5.6.1]{kahn_levine}.
Kahn and Levine show that this actually becomes a complex, $z_s(X,*;\mathcal{A})$, and they define
the higher Chow groups with coefficients in $\mathcal{A}$ as
\begin{equation*}
    \CH_s(X,t;\mathcal{A})=\Hoh_t(z_s(X,*;\mathcal{A})).
\end{equation*}

There are maps relating the complex $z_r(X,*;\mathcal{A})$ to $z_r(X,*)$, the untwisted
complex that computes Bloch's higher Chow groups. These are induced by the base-change map $B_\mathcal{E}$ and the forgetful map $F$ on $K$-theory.
\begin{align*}
    B_\mathcal{E}   &:\K_0(k(W))                \rightarrow    \K_0(k(W);\mathcal{A})\\
    F   &:\K_0(k(W),\mathcal{A})    \rightarrow    \K_0(k(W))
\end{align*}
The map $B_\mathcal{E}$ takes a $k(W)$-vector space and tensors with $\mathcal{E}_{k(W)}$ to produce a left $\mathcal{A}_{k(W)}$-module.
The norm map $F$ simply forgets the $\mathcal{A}\otimes_{k(W)}$-module structure on a vector space.

In particular, the kernels of these maps are zero, and the cokernels of the maps are
\begin{align}\label{eq:cokernels}
    coker(B_\mathcal{E})   &\iso   \ZZ/\left(\frac{na}{ind(\mathcal{A}_{k(W)})^2}\right)\\
    coker(F)   &\iso   \ZZ/\left(ind(\mathcal{A}_{k(W)})^2\right)
\end{align}
over $k(W)$.

\begin{lemma}\label{lem:multiplication}
    The compositions $F\circ B_\mathcal{E}$ and $B_\mathcal{E}\circ F$ are multiplication by $na$ on $z_s(X,t)$ and
    $z_s(X,t;\mathcal{A})$.
    \begin{proof}
        This follows immediately from Equation~\eqref{eq:cokernels}.
    \end{proof}
\end{lemma}

\begin{corollary}
    The cokernel of $F:z_s(X,t;\mathcal{A})\rightarrow z_s(X,t)$ is a torsion group of exponent bounded above by $n^2$,
    and $B_{\mathcal{E}}:z_s(X,t)\rightarrow z_s(X,t;\mathcal{A})$ is a torsion group of exponent bounded above by $na$.
    \begin{proof}
        In the first case, one always has $ind(\mathcal{A}_{k(W)})|n$, so the statement follows from Equation~\eqref{eq:cokernels}.
        Similarly,
        \begin{equation*}
            \left(\frac{na}{ind(\mathcal{A}_{k(W)})^2}\right)|na,
        \end{equation*}
        so the second statement follows.
    \end{proof}
\end{corollary}

\begin{proposition}
        The kernels and cokernels of
            \begin{equation*}
                B_\mathcal{E}:\CH_s(X,t)\rightarrow\CH_s(X,t;\mathcal{A})
            \end{equation*}
            and of
            \begin{equation*}
                F:\CH_s(X,t;\mathcal{A})\rightarrow\CH_s(X,t)
            \end{equation*}
            are torsion groups of exponent at most $na$.
    \begin{proof}
        This follows immediately from Lemma~\ref{lem:multiplication}.
    \end{proof}
\end{proposition}

Here is our main theorem.

\begin{theorem}
    Let $X$ be a $d$-dimensional scheme of finite type over a field, and let $\mathcal{A}$ be an Azumaya algebra on $X$.
    Then, the kernels and cokernels of
    \begin{equation*}
        B_\mathcal{E}:\G_r(X)\rightarrow \G_r(X;\mathcal{A})
    \end{equation*}
    and of
    \begin{equation*}
        F:\G_r(X;\mathcal{A})\rightarrow \G_r(X)
    \end{equation*}
    are groups of exponent bounded above by $(na)^{d+1}$ for all $r\geq 0$.
    \begin{proof}
        Kahn and Levine~\cite{kahn_levine} show that there is a convergent spectral sequence
        \begin{equation*}
            \Eoh_2^{p,q}(\mathcal{A})=\CH_q(X,-p-q;\mathcal{A})\Rightarrow \G_{-p-q}(X;\mathcal{A}).
        \end{equation*}
        There is also the motivic spectral sequence
        \begin{equation*}
            \Eoh_2^{p,q}=\CH_q(X,-p-q)\Rightarrow \G_{-p-q}(X).
        \end{equation*}
        The functors $B_\mathcal{E}:\G(X)\rightarrow \G(X;\mathcal{A})$ and $F:\G(X;\mathcal{A})\rightarrow \G(X)$
        are compatible with these spectral sequences and the functors $B_\mathcal{E}$ and $F$ on higher Chow
        groups. Note that $\Eoh_2^{p,q}=\Eoh_2^{p,q}(\mathcal{A})=0$ whenever $q<0$, $-p<0$, or $q>d$.

        We will prove the theorem for the kernel of the functor $B_\mathcal{E}$. The other cases are entirely similar. On the
        $\Eoh_{\infty}$-page, the composition functor $F\circ B_\mathcal{E}$ is still multiplication by $na$,
        so the kernels and cokernels of $B_\mathcal{E}$ on $\Eoh_{\infty}$ are still of exponent at most $na$.
        The spectral sequences abut to filtrations $F^s\G_r(X;\mathcal{A})$ and $F^s\G_r(X)$ where
        \begin{gather*}
            F^{(s/s+1)}\G_r(X;\mathcal{A})=F^s\G_r(X;\mathcal{A})/F^{s+1}\G_r(X;\mathcal{A})\iso\Eoh_{\infty}^{-r+s,-s}(\mathcal{A})\\
            F^{(s/s+1)}\G_r(X)=F^s\G_r(X)/F^{s+1}\G_r(X)\iso\Eoh_{\infty}^{-r+s,-s}.
        \end{gather*}
        The filtration looks like
        \begin{equation*}
            0=F^{0}\G_r(X)\subseteq F^{-1}\G_r(X)\subseteq\cdots\subseteq F^{-d}\G_r(X)=\G_r(X).
        \end{equation*}
        The filtration $F^s\G_r(X)$ is of length $d$ by the vanishing statements. Let $z\in \G_r(X)$
        be in the kernel of $F$, and let $\overline{z}$ be the image of $z$ in $E_{\infty}^{-r-d,d}$.
        Then, by hypothesis, $\overline{z}$ is in the kernel of $F$, so that $na\cdot\overline{z}=0$.
        Thus, $na\cdot z$ is contained in $F^{-d+1}\G_r(X)$.
        Continuing in this way, we see that $(na)^{d+1}\cdot z$ is contained in $F^0\G_r(X)=0$. So,
        $(na)^{d+1}\cdot z=0$.
    \end{proof}
\end{theorem}

\begin{corollary}
    The same result holds for $K$-theory when $X$ is regular.
\end{corollary}

\begin{corollary}
    The maps
    \begin{gather*}
        B_{\mathcal{E}}:F^{(s/s+1)}\G_r(X;\mathcal{A})\rightarrow F^{(s/s+1)}\G_r(X)\\
        F:F^{(s/s+1)}\G_r(X)\rightarrow F^{(s/s+1)}\G_r(X;\mathcal{A})
    \end{gather*}
    have torsion kernels and cokernels of exponent at most $na$.
    \begin{proof}
        This follows from the proof of the theorem.
    \end{proof}
\end{corollary}

\begin{corollary}
    For any integer $j$ prime to $na$, the maps
    \begin{gather*}
        B_\mathcal{E}:z_s(X,*;\ZZ/j)\rightarrow z_s(X,*;\mathcal{A};\ZZ/j)\\
        B_\mathcal{E}:\G_r(X;\ZZ/j)\rightarrow \G_r(X;\mathcal{A};\ZZ/j)\\
        F:z_s(X,*;\mathcal{A};\ZZ/j)\rightarrow z_s(X,*;\ZZ/j)\\
        F:\G_r(X;\mathcal{A};\ZZ/j)\rightarrow\G_r(X;\ZZ/j)
    \end{gather*}
    are isomorphisms.
\end{corollary}

It is interesting that this method proves the isomorphisms by means of an isomorphism of cycle complexes, not just a quasi-isomorphism.

\bibliographystyle{amsplain-pdflatex}
\bibliography{brauer,../bibliographies/mypapers}

\noindent
Benjamin Antieau
[\texttt{\href{mailto:antieau@math.ucla.edu}{antieau@math.ucla.edu}}]\\
UCLA\\
Math Department\\
520 Portola Plaza\\
Los Angeles, CA 90095-1555\\

\end{document}